\documentclass{amsart}
\usepackage{amsmath,amsthm}
\usepackage{amsfonts,amssymb}

\newtheorem{thm}{Theorem}[section]
\newtheorem{cor}[thm]{Corollary}
\newtheorem{lem}[thm]{Lemma}
\newtheorem{prop}[thm]{Proposition}

\newtheorem{defn}[thm]{Definition}

\theoremstyle{remark}

 \def\CO{{\mathcal O}}

 \def\RR{{\mathbb R}}

\begin{document}
 
\title[mean convergence of interpolating polynomials]
{Mean convergence of orthogonal Fourier series and interpolating polynomials}
\author{P\'eter V\'ertesi}
\address{R\'enyi Mathematical Institute\\ Hungarian Academy of Science\\
     Re\'altanoda U. 13-14\\1053 Budapest\\
    Hungary}\email{veter@renyi.hu}
\author{Yuan Xu}
\address{Department of Mathematics\\ University of Oregon\\
    Eugene, Oregon 97403-1222.}\email{yuan@math.uoregon.edu}

\date{August 10, 2003}
\thanks{The work of the first author was supported by Hungarian National 
Foundation for Scientific Research No. T22943, T32872 and T37299; the work of 
the second author was supported in part by the National Science Foundation 
under Grant DMS-0201669}
                                                    
\begin{abstract}
For a family of weight functions that include the general Jacobi weight 
functions as special cases, exact condition for the convergence of the Fourier
orthogonal series in the weighted $L^p$ space is given. The result is then 
used to establish a Marcinkiewicz-Zygmund type inequality and to study weighted
mean convergence of various interpolating polynomials based on the zeros of 
the corresponding orthogonal polynomials.
\end{abstract}

\maketitle                      
 
\section{Introduction} 
\setcounter{equation}{0}

Let $d \alpha$ be a finite nonnegative measure on $[-1,1]$. We consider the
Fourier orthogonal expansion with respect to $d\alpha$ and weighted $L^p$ 
convergence of the interpolation polynomials based on the zeros of orthogonal
polynomials with respect to $d\alpha$. 

Throughout this paper we denote by $L^p(d\alpha)$ the space of measurable
function $f$ such that 
$$
\|f\|_{d\alpha,p} = \left(\int_{-1}^1 |f(x)|^p d\alpha(x) \right)^{1/p}, 
\qquad    0 < p < \infty, 
$$
is finite. We assume that $\|f\|_\infty$ is the usual uniform norm for the 
continuous functions. If $d\alpha = w dt$, we may write $\|f\|_{w,p}$ instead 
of $\|f\|_{d\alpha,p}$. Let $p_n(d\alpha)$ denote the orthonormal polynomial 
of degree $n$ with respect to the measure $d\alpha$ on $[-1,1]$. The zeros 
of $p_n(d\alpha)$ are distinct real numbers, denoted by $x_{1n}(d\alpha),
x_{2n}(d\alpha),\ldots, x_{nn}(d\alpha)$, in $(-1,1)$. For any given function
$f$ on $[-1,1]$, we let $L_n(d\alpha;f)$ denote the unique Lagrange 
interpolation polynomial of degree $n-1$ that agrees with $f$ at 
$x_{kn}(d\alpha)$, $1 \le k \le n$.
 
Let $d\beta$ be another measure defined on $[-1,1]$. We are interested in 
the precise condition on $d\beta$ and $d\alpha$ that will ensure the 
convergence of $L_n(d\alpha;f)$ in the $L^p(d\beta)$ norm for $f\in C[-1,1]$.
This question was addressed by many authors (see \cite{MV,N,Ve,VeXu,X} for 
historical account). In \cite{N}, Nevai solved the problem for the case that 
$\alpha'$ and $\beta'$ are generalized Jacobi weight functions, defined by 
\begin{equation} \label{eq:1.1}
w(x) = h(x) \prod_{i=0}^{r+1} |x-t_i|^{\Gamma_i}, \qquad
 -1 = t_0 < t_1 < \ldots < t_r < t_{r+1} = 1,
\end{equation}
where $h$ is a positive continuous function on $[-1,1]$ and the modulus of 
continuity $\omega$ of $h$ satisfies $\int_0^1(\omega(t)/t) dt < + \infty$. 
The condition for the convergence of $L_n(d\alpha;f)$ in $L^p(d\beta)$, 
$0 < p < \infty$, is given by 
\begin{equation} \label{eq:1.2}
(\alpha' \sqrt{1-x^2})^{-p/2} \beta' \in L^1.
\end{equation} 
Since then this result has been extended in several directions, to other
interpolation process and to more general weight functions (see, for example, 
\cite{MV,N,Ve,VeXu,X} and the reference therein). It turned out (\cite{X}) 
that one way of proving such results is to use a Marcinkiewicz-Zygmund type 
inequality. In the simplest case, such an inequality takes the form of 
$$
\| P\|_{d\beta,p} \le c \left(\sum_{k=1}^n c_{k,n}
  |P\left(x_{kn} (d\alpha)\right)|^p \right)^{1/p}, 
$$
where $P$ is any polynomial of degree at most $n-1$ and $c_{k,n}$ are certain
(precisely known) nonnegative numbers. For $\alpha'$ and $\beta'$ being 
generalized Jacobi weight functions, the inequality was proved in \cite{X} 
under the precise condition \eqref{eq:1.2} for the convergence of 
$L_n(d\alpha;f)$; furthermore, it was extended to include derivative values 
in the right hand side so that the convergence of Hermite interpolation 
polynomials can be derived.

We will try to establish the Marcinkiewicz-Zygmund type inequality for more
general weight functions. To be sure, this has been done in \cite{MV}; but
the result there requires an additional condition other than \eqref{eq:1.2},
namely, $(\alpha'\sqrt{1-x^2})^{q/2} \beta'^{q-1} \in L^1$, where $p^{-1}+
q{-1} =1$. To establish the inequality, one way is to use the boundedness of 
the orthogonal Fourier expansion. For $f \in L^2(d\alpha)$, such an expansion 
is given by
\begin{equation} \label{eq:1.3}
 f \sim \sum_{n=0}^\infty c_n(f) p_n(d\alpha), \qquad c_n(f) = 
    \int_{-1}^1 f(t) p_n(d\alpha;t) d\alpha.
\end{equation}
Let $S_n(d\alpha;f)$ denote the $n$-th partial sum of the expansion. The
convergence of the Fourier expansion amounts to the uniform boundedness of 
$S_n(d\alpha;f)$. Finding the precise conditions for the boundedness of 
$\|S_n(d\alpha;f)\|_{d\beta,p}$ is an interesting problem in itself. By the 
Christoeff-Darboux formula, the kernel $K_n(x,y)$ of this operator has a 
singularity at $x=y$. To overcome the problem of the singularity, \cite{X} 
established the following inequality 
$$
\int_{-1}^1 \Big| \int_{-1}^1 \frac{g(y)}{x-y}\, dy\Big|^p
U^p(x) \, dx \le c\int_{-1}^1 |g(x)|^p V^p(x)\, dx,
$$
in which $U$ and $V$ are the generalized Jacobi weight functions in 
\eqref{eq:1.1}. This is an inequality for the Hilbert transform defined by
\begin{equation*} 
  H(g;x)= \lim_{\varepsilon \to 0} \int_{|x-y|\ge \varepsilon} 
    \frac{g(y)}{x-y} dy, \qquad g \in L^1.
\end{equation*}
In \cite{MW}, double weight inequality for the Hilbert transform on 
$[0,\infty)$ is proved for general weight functions, the so-called $A_p$ 
weight. Although the results in \cite{MW} do not apply to weight functions 
\eqref{eq:1.1} directly, we show that they can be used to establish 
inequalities for weight functions that are more general than the weight 
function in \eqref{eq:1.1} (see, for example, (2.2) and Definition 3.6),
which in turn gives the result on mean convergence of various interpolating
polynomials. 

The paper is organized as follows: In the next section we fix the notation
and state the preliminary. In section 3 we prove a double weight inequality 
for the Hilbert transform. The Fourier orthogonal series is studied in 
Section 4. The Marcinkiewicz-Zygmund inequalities are proved in
Section 5, followed by the discussion on the mean convergence of interpolating
polynomials in Section 6.

\section{Notation and Preliminary} 
\setcounter{equation}{0}

Throughout this paper we denote by $\Pi_n$ the space of polynomials of 
degree at most $n$ and by $\Pi$ the space of all polynomials. We will
use constants $c,c_1,c_2\ldots$ to denote generic constants that depend 
only on weight functions and other fixed parameters involved, their values
may vary from line to line. The notation $A\sim B$ means $|A^{-1}B|\le c$ 
and $|AB^{-1}| \le c$.

\subsection{Weight function} First we define the weight functions that we 
shall deal with in this paper.

\begin{defn} \label{def2.1}
A function $w$ is called a generalized Jacobi weight function ($w \in GJ$ in 
short) if, for $t \in (-1,1)$, 
\begin{equation} \label{eq:2.1}
 w(t) = h(t) \prod_{i=0}^{r+1} \left[\tau_i(|t-t_i|)\right]^{a_i},
 \qquad  - 1 = t_0 < t_1 < \ldots < t_r < t_{r+1} =1,
\end{equation}
where $a_i \in \RR$ and $\tau_i$ are nondecreasing, continuous semi-additive
functions, $\tau_i(0)=0$, $h$ is a nonnegative function that satisfies 
$h \in L^\infty[-1,1]$ and $1/h \in L^\infty[-1,1]$ (we do not assume that
$w \in L^1[-1,1]$). 

A measure $d\alpha = \alpha' dt$ is called a GJ measure if 
$\alpha' \in GJ$ and $\alpha' \in L^1[-1,1]$. 
\end{defn}

Throughout the paper we often write $w \in GJ$ as 
\begin{equation} \label{eq:2.3}
 w(t) = h(t)\prod_{i=0}^{r+1} w_i(|t-t_i|), 
    \qquad w_i(t) = \left[\tau_i(t)\right]^{a_i},
\end{equation}
and impose conditions on $w_i(t)$ instead on $\tau_i$.

\begin{defn} \label{def2.2}
Let $w \in GJ$ as in Definition \ref{def2.1}. For $i = 0,1,\ldots, r+1$, 
consider
\begin{enumerate}
\item $\tau_i(t)$ is concave; that is, $\tau_i(t)+\tau_i(s) \le 
 2\tau_i((t+s)/2)$; 
\item
$\displaystyle{\int_0^\delta w_i(s)ds = \CO(\delta w_i(\delta))}, 
\quad \delta \to +0$;
\item $ \omega(h;\delta)_\infty \delta^{-1} \in L^1[0,1] \quad \hbox{or} \quad
  \omega(h,\delta)_2 = \CO(\sqrt{\delta}), \quad \delta \to +0$, where $\omega$
  is a modulus of continuity.
\end{enumerate}
We say $w \in GJ1$ if it satisfies (1), $w \in GJ2$ if it satisfies (2), and
$w \in GJ3$ if it satisfies (3). We call $w$ an admissible GJ weight function
if it satisfies all three conditions.
\end{defn}

If all $\tau_i(t) = t$, then $w \in GJ$ is precisely the usual generalized 
Jacobi weight function in \eqref{eq:1.1}, which is an admissible GJ weight 
function. Another family of GJ weight functions is given as follows.

\begin{defn}\label{def:log}
Let $\Gamma_i$ and $\gamma_i$ be real numbers. We denote by $GJ\log$ the 
collection of the GJ weight functions
\begin{equation} \label{eq:2.2}
  w(t) = h(t) \prod_{i=0}^{r+1} |t-t_i|^{\Gamma_i} \log^{\gamma_i} 
      \frac{e}{|t-t_i|}. 
\end{equation} 
\end{defn}

The weight function $w \in GJ\log$ is in GJ2 if $\Gamma_i > -1$, $0\le i \le 
r+1$ (\cite[p. 328]{MV}). Furthermore, $\Gamma_i > -1$, $0\le i \le r+1$ 
is necessary for $w \in GJ2$ and $w \in L^1$. We note that this condition is
slightly stronger than $w \in L^1$. Also, $w \in GJ\log$ is in GJ1 if, say, 
$\gamma_i/\Gamma_i \ge 0$. We often write $\Gamma_i(w)$ or $\Gamma_i(d\alpha)$ 
in place of $\Gamma_i$ to emphasis that they are parameters of $w$ or 
$d\alpha$, respectively, of the above form. 

For some results we also need one more restriction on the weight functions. 

\begin{defn} \label{def:2.3}
If $w \in GJ$ can be written as $w = u/v$ such that both $u$ and $v$ are 
in the form of \eqref{eq:2.1} with positive exponents, and $v$ satisfies 
$$
\int_0^\delta \frac{1}{v_i(s)} ds = 
  \CO\left(\frac{\delta}{v_i(\delta)}\right),
 \qquad 0 \le i \le r+1, \quad
     v(t) = h(t) \prod_{i=0}^{r+1} v_i(|t-t_i|),
$$
then we say that $w \in GJ4$.
\end{defn}

As an example, we point out that weight functions $w$ in \eqref{eq:2.2} are 
GJ4 weight functions if $\Gamma_i > -1$.

Throughout this paper, we reserve the notation $\varphi$ for the functions
$$
 \varphi(x) = \sqrt{1-x^2} \quad \hbox{and} \quad 
        \varphi(n,x)=\sqrt{1-x^2}+n^{-1}, \quad
 -1 \le x \le 1.
$$
Furthermore, for $w \in GJ$ as in \eqref{eq:2.3}, we define 
\begin{equation} \label{eq:2.4}
 w(n,t) =\frac{w_0(\sqrt{1-t}+n^{-1})w_{r+1}(\sqrt{1+t}+n^{-1})}
 {\sqrt{1-t^2}+ n^{-1}} \prod_{i=1}^r w_i(|t-t_i|+n^{-1}).
\end{equation}

\subsection{Orthogonal polynomials}
We consider orthonormal polynomials $p_n(d\alpha) \in \Pi_n$ with respect to 
$d\alpha = \alpha' dt$ and $\alpha' \in GJ$; that is, 
$$
  \int_{-1}^1 p_n(d\alpha;t)p_m(d\alpha;t) d\alpha = \delta_{n,m}. 
$$
We assume that the zeros $x_{kn}(d\alpha)$ of $p_n(d\alpha)$ take the order 
$$
  -1 < x_{nn}(d\alpha) < x_{n-1,n}(d\alpha) < \ldots < x_{1,n}(d\alpha) < 1.
$$ 
The Christoeff function $\lambda_n(d\alpha)$ with respect to $d\alpha$ is 
defined by 
$$
  \lambda_n(d\alpha;t) = \min_{P \in \Pi_{n-1}} \frac{1}{|P_n(t)|^2}
    \int_{-1}^1 |P(x)|^2 w(x)dx. 
$$
The numbers $\lambda_{kn}(d\alpha)=\lambda_n(d\alpha; x_{kn})$ are called the 
Cotes numbers, which appear in the Gauss quadrature formula
\begin{equation*}
\sum_{k=1}^n P\bigl(x_{kn}(d\alpha)\bigr) \lambda_{kn}(d\alpha)
   = \int_{-1}^1 P(t) d\alpha, \qquad p \in \Pi_{2n-1}. 
\end{equation*}
Further assumption on the weight function is needed to get bounds on these
quantities. 

\begin{lem} \label{lem:2.4}
Let $d\alpha$ be an admissible GJ measure. Then
$$
|p_n (d\alpha;x)| \le c\,\alpha'(n,x)^{-1/2} \varphi(n,x)^{-1/2}
$$
uniformly for $-1\le x\le 1$ and 
$$
|p'_n (d\alpha; x_{kn})| \sim \frac{n}{\varphi(n,x_{kn})} \frac{1}{
   \bigl[\alpha'(n,x_{kn}) \varphi(n,x_{kn})\bigr]^{1/2} }
$$
uniformly for  $1\le k\le n$, where $x_{kn} = x_{kn} (d\alpha)$, and
$$
\lambda_{n} (d\alpha;x) \sim \frac{1}{n} \alpha'(n,x) \varphi(n,x)
$$
uniformly for $-1 \le x \le 1$. 
\end{lem}

For the proof of these estimates, see \cite{Ba,MV}. It should be pointed out
that some of the estimates hold for more general weight functions, or hold
under weaker conditions for the weight functions in Definition \ref{def2.2}.
For example, for the estimate of $\lambda_n(d\alpha;x)$, only $\alpha' \in
GJ2$ is needed. See the discussions in \cite{MV}.

Let $w\in GJ$, for a fixed $d>0$, we define $\Delta_n(\varepsilon)$ by
\begin{equation} \label{eq:2.5}
\Delta_n (\varepsilon) = [-1 + \varepsilon n^{-2},1 - \varepsilon n^{-2}] 
 \Big \backslash 
 \bigcup_{i=1}^r  [t_i- \varepsilon n^{-1},\, t_i+ \varepsilon n^{-1}].  
\end{equation}
We shall use $\chi_E$ to denote the characteristic function of a set $E$.
The following lemma is a simplified version of Theorem 3.5 in \cite{MV}.  

\begin{lem} \label{lem:2.5}
Let $d \beta$ be a GJ2 measure and $u \in GJ$. Then for each $0<p<+\infty$ 
there exists an $\varepsilon_0 >0$ such that for every fixed $\varepsilon$,
$0 < \varepsilon \le \varepsilon_0$, and for $P\in \Pi_n$,
$$
\int_{-1}^1 |P(t)|^p u(n,t) d\beta \le c 
   \int_{\Delta_n(\varepsilon)} |P(t)|^p u(n,t) d\beta, \qquad n \ge n_0.  
$$
\end{lem}

The next lemma gives an inequality for the quadrature sum of polynomials:

\begin{lem} \label{lem:2.6}
Let $\alpha'$ be an admissible GJ measure and $v \in GJ2 \cap GJ4$. Then for 
$1 \le p < \infty$, 
$$
  \sum_{k=1}^n \lambda_{kn}\left(v,x_{kn}(d\alpha)\right)
  \left|P(x_{k,n}(d\alpha))\right|^p 
   \le c \int_{-1}^1 |P(x)|^p v(x) dx 
$$
for every $P \in \Pi_{mn}$, where $m$ is a fixed positive integer and $c$ is
independent of $P$ and $n$.
\end{lem}

Finally, there is the inequality of Bernstein-Markov type for general weight
functions.

\begin{lem}\label{lem:2.7}
Let $d\alpha$ be a GJ4 measure and $w \in GJ$. Let $0 < p < \infty$. Then
for arbitrary $P \in \Pi_n$ and integer $j \ge 1$,
\begin{align*}
  \int_{-1}^1 \left|P^{(j)}(x) \varphi^j(n,x)\right|^p w(n,x) \varphi(n,x) 
   d\alpha \le c n^{jp} \int_{-1}^1 |P(x)|^p w(n,x) \varphi(n,x)d\alpha.
\end{align*} 
\end{lem}

The last two results have been studied by several authors for various 
weight functions. In its present generality, it appears as 
\cite[Theorem 3.6]{MV} and \cite[Theorem 2.D]{MV}. See also \cite{MT}
for some of the above inequalities with doubling weight.

\section{Weighted inequalities for the Hilbert transform} 
\setcounter{equation}{0}

We start with a result on the Hilbert transform proved in \cite{MW} for 
general weight functions defined on $[0, \infty)$. 

\begin{lem}\label{lem:3.1}
Let $U$ and $V$ be nonnegative weight functions defined on $[0, \infty)$,
and there exists a constant $A$ such that either
\begin{equation} \label{eq:3.1}
U(x) \le A U(y) \quad \hbox{and} \quad  V(x) \le A V(y), \qquad 
  x \le y \le 2x, \quad x >0, 
\end{equation}
or the similar inequalities with $\ge$ in place of $\le$ hold.
Let $1 < p < \infty$ and $f V \in L^p$. Then there is a constant $c$ 
independent of $f$ such that 
\begin{equation} \label{eq:3.2}
\int_{0}^\infty \left| \int_0^\infty \frac{f(y)}{x-y} dy \right|^pU^p(x) dx
   \le c \int_{0}^\infty |f(x)|^p V^p(x) dx, 
\end{equation} 
if for every interval $I \subset [0,\infty)$,  
\begin{equation} \label{eq:3.3}
 \left[ \int_{0}^\infty  \frac{|I|^{p-1} [U(x)]^p dx}{(|I|+|x-x_I|)^p} \right]
 \left[\frac{1}{|I|} \int_I [V(x)]^{-q} dx \right]^{p-1} \le B,
\end{equation} 
and 
\begin{equation} \label{eq:3.4}
\left[\frac{1}{|I|} \int_I [U(x)]^p dx \right]  
\left[\int_{0}^\infty \frac{|I|^{q-1}[V(x)]^{-q}}{(|I|+|x-x_I|)^{q}} 
   dx \right]^{p-1} \le B,
\end{equation} 
where $q= p/(p-1)$, $|I|$ denotes the length of $I$, $x_I$ is the center 
of $I$ and $B$ is independent of $I$.
\end{lem}

The lemma is stated in \cite[Theorem 8]{MW} with the integral of the Hilbert 
transform in the left hand side. The proof there shows that the above version 
holds. By translation and truncation, it is possible to state a version of this
theorem for the interval $[-1,1]$. However, the condition \eqref{eq:3.1} 
does not hold for GJ weight function. Our main result in this section is 
a double weight inequality for the GJ weight functions. The following two
lemmas will be useful (\cite[p. 281 and 282]{MW}).

\begin{lem}\label{lem:3.2}
Let $1 \le p < \infty$ and $fV \in L^1$. There is a finite $c$, independent of 
$f$, such that 
$$
\int_{0}^\infty \left|\int_0^x f(t)dt \right|^p U^p(x) dx 
   \le c \int_{0}^\infty |f(x)|^p V^p(x) dx 
$$
if and only if there is a finite $B$, independent of $\delta$, such that for
$\delta > 0$, 
\begin{equation} \label{eq:3.5}
\left[\int_\delta^\infty [U(x)]^p dx \right]  
\left[\int_{0}^\delta [V(x)]^{-q} dx \right]^{p-1} \le B.
\end{equation} 
\end{lem}

\begin{lem}\label{lem:3.3}
Let $1 \le p < \infty$ and $fV\in L^p$. There is a finite $c$, independent of 
$f$, such that 
$$
\int_{0}^\infty \left|\int_x^\infty f(t)dt \right|^p U^p(x) dx 
   \le c \int_{0}^\infty |f(x)|^p V^p(x) dx 
$$
if and only if there is a finite $B$, independent of $\delta$, such that for
$\delta > 0$, 
\begin{equation} \label{eq:3.6}
\left[\int_0^\delta [U(x)]^p dx \right]  
\left[\int_{\delta}^\infty [V(x)]^{-q} dx \right]^{p-1} \le B.
\end{equation} 
\end{lem}

Our main result in this section is the following theorem, in which $U_i$ are
parts of $U$ as in the notation \eqref{eq:2.3}.

\begin{thm}\label{thm:3.4}
Let $U$ and $V$ be GJ weight functions. Let $1 < p < \infty$ and $g V \in L^p$.
Then there is a constant $c$ independent of $f$ such that 
\begin{equation} \label{eq:3.7}
\int_{-1}^1 \left| \int_{-1}^1 \frac{g(y)}{x-y} dy \right|^p U^p(x) dx
  \le c \int_{-1}^1 |g(x)|^p V^p(x) dx, 
\end{equation} 
if there is a $B$, independent of $\delta$, such that for $\delta > 0$,
\begin{equation} \label{eq:3.8}
 \left[ \int_{0}^1  \frac{U_i^p(t) dt}{(\delta + t)^p} \right]
 \left[\int_0^\delta V_i^{-q}(t)dt \right]^{p-1} \le B, 
\qquad 0 \le i \le r+1, 
\end{equation} 
and 
\begin{equation} \label{eq:3.9}
\left[ \int_0^\delta U_i^p(t) dt \right]  
\left[\int_{0}^1 \frac{V_i^{-q}(t)}{(t +\delta)^q} dt \right]^{p-1} \le B,
\qquad 0 \le i \le r+1.
\end{equation} 
\end{thm}

\begin{proof}
Recall that $t_i$ are fixed numbers, $-1=t_0 < t_1 < \ldots < t_r <t_{r+1} =1$.
We write
\begin{align*}
&\int_{-1}^1 \Big| \int_{-1}^1 \frac{g(y)}{x-y}\, dy\Big|^p
U^p(x) \, dx  = \sum_{i=0}^r \int_{t_i}^{t_{i+1}} \Big| \int_{-1}^1
\frac{g(y)}{x-y}\, dy\Big|^p U^p(x)\, dx \cr
&\qquad \le c\sum_{i=0}^r \int_{t_i}^{t_{i+1}} \Big| 
  \int_{-1}^1 \frac{g(y)}{x-y}\,dy\Big|^p U_i^p(|x-t_i|) 
  U_{i+1}^p(|x-t_{i+1}|) dx.
\end{align*}
For each $i$, we then break the inner integral of the last expression into
three integrals over $(-1,t_i)$, $(t_i,t_{i+1})$ and $(t_{i+1},1)$, 
respectively, and estimate the corresponding terms separately. 

We estimate the middle part first. Changing variables
$$
x-t_i = \frac{t_{i+1}-t_i}{1+X}\quad ,\quad y-t_i = \frac{t_{i+1} - t_i}
{1+Y}
$$
in the integrals gives 
\begin{align*}
& J_i:= \int_{t_i}^{t_{i+1}} \Big| \int_{t_i}^{t_{i+1}}
\frac{g(y)}{x-y}\, dy\Big|^p U_i^p(|x-t_i|) U_{i+1}^p(|x-t_{i+1}|) dx \\
& = \int_0^{+\infty} \Big| \int_0^{+\infty} \frac{g(y)}{X-Y} 
 \frac{X+1}{Y+1}  dY\Big|^p U_i^p\left(\frac{t_{i+1}-t_i}{X+1}\right)
    U_{i+1}^p\left(\frac{t_{i+1}-t_i}{X+1} X\right)
   \frac{{t_{i+1} -t_i}}{(X+1)^2} dX.
\end{align*}
If we can apply Lemma \ref{lem:3.1} with $f(Y)= g(y)/(1+Y)$ and with $u$ 
and $v$ in place of $U$ and $V$, where  
\begin{align*}
& u^p(X): = U_i^p\left(\frac{t_{i+1}-t_i}{X+1}\right) 
  U_{i+1}^p\left(\frac{t_{i+1}-t_i}{X+1} X\right)
  (X+1)^{p-2} \\ 
& v^p(X): = V_i^p\left(\frac{t_{i+1}-t_i}{X+1} \right)
 V_{i+1}^p\left(\frac{t_{i+1}-t_i}{X+1} X\right)(X+1)^{p-2},  
\end{align*}
then we will end up with the desired estimate
\begin{align*}
J_i & \le c \int_0^{+\infty} |g(x)|^p\Big|^p 
  V_i^p\left(\frac{t_{i+1}-t_i}{X+1}\right)
  V_{i+1}^p\left(\frac{t_{i+1}-t_i}{X+1} X\right) 
  \frac{{t_{i+1} -t_i}}{(X+1)^2} dX \\
& = c \int_{t_i}^{t_{i+1}} |g(y)|^p V_i^p(|y-t_i|) V_{i+1}^p(|y-t_{i+1}|)dy
  \le c \int_{t_i}^{t_{i+1}} |g(y)|^p V^p(y)dy. 
\end{align*}
The functions $U_i$ and $V_i$ satisfy \eqref{eq:3.1} automatically, since
$\tau_i$ in the Definition \ref{def2.1} is nondecreasing so that $U_i$ and
$V_i$ are either nondecreasing or nonincreasing. Thus, we only need to 
verify that the conditions \eqref{eq:3.3} and \eqref{eq:3.4} of Lemma 
\ref{lem:3.1} are satisfied with $u$ and $v$ in place of $U$ and $V$. First 
we note that the change of variables leads to $u^p(x)/(1+X)^p \in
L^1[0,\infty)$ and $v^{-q}(x)/(1+X)^q \in L^1[0,\infty)$, since \eqref{eq:3.8} 
and \eqref{eq:3.9} with a fixed $\delta$ shows that $U^p \in L^1$ and 
$V^{-q}\in L^1$. Hence, if $|I|$ is fixed, then the inequality $|I| + 
|X - X_I| \ge c (1+X)$ shows that the \eqref{eq:3.3} and \eqref{eq:3.4} 
holds trivially. Thus, we need to consider only the following two cases:

\medskip
\noindent
{\it Case 1.} $I = (0, 2 \delta)$, $\delta < 1/2$. Using the fact that
$|I| + |X- X_I| = 2 \delta + |X- \delta| \ge (1+X)/2$ if $X >1$, the
left hand side of \eqref{eq:3.3} is bounded by
\begin{align*}
&  c + \left(\int_0^1 \frac{[u(X)]^p}{(\delta + |X- \delta|)^p} dX\right)
   \left(\int_0^{2\delta} [v(X)]^{-q} dX \right)^{p-1}\\
\le \; & c + c \left(\int_0^1 \frac{\left[U_{i+1}\left(\frac{t_{i+1}-t_i}{1+X} 
 X\right)\right]^p} {(2\delta + |X- \delta|)^p} dX\right)
   \left(\int_0^{2\delta} \left[V_{i+1}\left(\frac{t_{i+1}-t_i}{1+X}X\right)
\right]^{-q} dX\right)^{p-1}\\
\le & \; c + c \left(\int_0^1\frac{\left[U_{i+1}(x)\right]^p} {(x+\delta)^p} dx
 \right)\left(\int_0^{2\delta} \left[V_{i+1}(x)\right]^{-q} dx\right)^{p-1}.
\end{align*}
The last expression is bounded by \eqref{eq:3.9}. The inequality 
\eqref{eq:3.4} in the case of $I \subset (0, 2\delta)$ is established 
similarly.

\medskip\noindent
{\it Case 2.} $I = (s, R)$, $s < R$, $R$ can be arbitrarily large. 
If $X \le R/2$, then $|I| + |X-X_I| = (R-s) + |X - (R+s)/2| \ge R-s/2\ge X+1$;
so that 
\begin{align*}
 \int_0^{R/2} \frac{[u(X)]^p}{(\delta + |X- \delta|)^p} dX
 \le \int_0^{R/2} \frac{[u(X)]^p}{(1+X)^p} dX \le c, 
\end{align*}
as $u^p/(1+X)^p \in L^1[0,\infty)$. The left had side of \eqref{eq:3.3} is 
bounded since changing variables from $X$ and $Y$ back to $x$ and $y$ shows 
that
\begin{align*}
  \int_{R/2}^\infty \frac{[u(X)]^p}{(\delta + |X- \delta|)^p} dX 
 & \le c \int_{R/2}^\infty \frac{[U_i\left(\frac{t_{i+1}-t_i}{1+X}\right)]^p}
     {(\delta + |X- \delta|)^p} (1+X)^{p-2} dX \\  
 & \le c \delta^p \int_{t_i}^{t_i+\delta} \frac{[U_i(|x-t_i)]^p}
        {(\delta + |x- t_i|)^p} dx 
   = c \int_{0}^{\delta}[U_i(t)]^p dt,      
\end{align*}
where $\delta = (t_{i+1}-t_i)/(1+R/2)$ and we have used the fact that
$$
(R-s) + |X - (R+s)/2| \ge R/2 + X-s \ge (R + X+1)/4,
$$
and since by $-q(p-2)/p=q-2$, 
\begin{align*}
& \int_s^R [v(X)]^{-q} dX  \le \int_s^R 
 \left[V_i\left(\frac{t_{i+1}-t_i}{1+X}\right)\right]^{-q} (1+X)^{q-2} dX \\  
 & \le c \int_{t_i + \delta}^{t_i + \rho} \left[V_i(|x-t_i|)\right]^{-q}
    \frac{dx}{|x-t_i|^q} \le c \int_{\delta}^{\rho} \left[V_i(t)\right]^{-q}
    \frac{dt}{t^q} \le c  \int_{0}^{1} \frac{\left[V_i(t)\right]^{-q}}
      {(t+\delta)^q} dt,
\end{align*}
where $\rho = (t_{i+1}-t_i)/(1+R)$ and $\delta$ is as above, so that the
boundedness of the left hand side of \eqref{eq:3.3} follows from 
\eqref{eq:3.9}. The inequality \eqref{eq:3.4} is established similarly.
Consequently, we have justified the use of Lemma \ref{lem:3.1} and the
bound of $J_i$. 

Next, let $\bar t_i = (t_i + t_{i+1})/2$. We split the first integral as 
follows,
\begin{align*}
&\int_{t_i}^{t_{i+1}} \left| \int_{-1}^{t_i} \frac{g(y)}{x-y} dy \right|^p 
  U_i(|x-t_i|) U_{i+1}^p(|x-t_{i+1}|)dx \\
&\le c\int_{\bar t_i}^{t_{i+1}} \left[\int_{-1}^{t_i} |g(y)| dy
 \right]^p U_{i+1}^p(|x-t_{i+1}|)dx\\
&\;\; + c\int_{t_i}^{\bar t_i} \left[ \int_{-1}^{\bar t_{i-1}}
|g(y)| dy\right]^p U_i^p(|x-t_i|) dx 
 + c \int_{t_i}^{\bar t_i} \left|\int_{\bar t_{i-1}}^{t_i}
\frac{g(y)}{x-y} dy \right|^p U_i^p(|x-t_i|) dx.
\end{align*}
The first two terms can be estimated by the H\"older inequality. For example,
the first term is bounded by
\begin{align*}
& c\int_{-1}^{t_i} |g(y)|^p V^p (y) dy \left[ \int_{-1}^{t_i}
\left[ V(y)\right]^{-q}\, dy\right]^{p/q} \int_{\bar t_i}^{t_{i+1}} 
  U_{i+1}^p(|x-t_{i+1}|) dx\\
&\qquad \le c \int_{-1}^{t_i} |g(y)|^p V^p(y) dy.
\end{align*}
The third term needs more work. Changing variables $x-t_i =X$, $y-t_i = -Y$, 
it becomes a constant multiple of 
$$
L:= \int_0^A \Big| \int_0^B \frac{g(y)}{X+Y} dY\Big|^p U_i^p(X) dX,
\qquad A = (t_{i+1} - t_i) /2, \quad B= (t_i-t_{i-1})/2.
$$
To estimate this term we need to use Lemma \ref{lem:3.2} and Lemma 
\ref{lem:3.3}. Changing variables 
$$
   X = \frac{A s}{1+s} \qquad \hbox{and} \qquad  Y = \frac{A t}{1+t}, 
$$
and splitting the inner integral into two parts gives 
\begin{align*}
 L & = A B^p \int_0^{\infty} \left| \int_0^{\infty}
 \frac{g(y)}{As(1+t) + Bt(1+s)}\frac{1+s}{1+t} dt\right|^p
 U_i^p\left(\frac{A s}{1+s}\right) \frac{ds}{(1+s)^2}\\ 
  & \le A^{1-p} B^p \int_0^{\infty} \left[ \int_0^s \frac{|g(y)|}
 {s(1+t)} \frac{1+s}{1+t} dt\right]^p U_i^p\left(\frac{As}{1+s}\right)
   \frac{ds}{(1+s)^2} \\
& \qquad + A \int_0^{\infty} \left[ \int_s^{\infty}
  \frac{|g(y)|}{t(1+s)}\frac{1+s}{1+t}dt\right]^p 
   U_i^p\left(\frac{A s}{1+s}\right) \frac{ds}{(1+s)^2}.
\end{align*}
To estimate the first integral in the right hand side, we use Lemma 
\ref{lem:3.2} with $f(t) = g(y) /(1+t)^2$ and with $u$ and $v$ in place of 
$U$ and $V$, where  
$$
  u^p(s) = U_i^p\left(\frac{As}{1+s}\right) \left(\frac{1+s}{s}\right)^p 
    \frac{1}{(1+s)^2} \quad \hbox{and}\quad 
   v^p(s) = V_i^p\left(\frac{Bt}{1+t}\right) (1+t)^{2(p-1)}, 
$$
so that the term is bounded by, after changing the integral back to $y$,
\begin{align*}
& c \int_0^\infty \left|\frac{g(y)}{(1+t)^2}\right|^p 
  V_i^p\left(\frac{Bt}{1+t}\right) (1+t)^{2p} \frac{dt}{(1+t)^2} \\
&   = c \int_0^B |g(y)|^p V_i^p(Y)dY
   = c \int_{\bar t_{i-1}}^{t_i} |g(y)|^p V_i^p(|y-t_i|)dy.
\end{align*}
The condition of Lemma \ref{lem:3.2} is verified as follows: changing 
variables back to $X$ and $Y$, 
\begin{align*}
 \int_\delta^\infty u^p(s) ds \left(\int_0^\delta v^{-q}(t) dt\right)^{p-1} 
&   = \left(\frac{A}{B}\right)^{p-1} \int_{A\delta}^A \frac{U_i^p(X)}{X^p} dX
     \left(\int_0^{B\delta} V_i^{-q}(Y)dY \right)^{p-1} \\
&   \le c \left(\frac{A}{B}\right)^{p-1} \int_{0}^A 
    \frac{U_i^p(x)}{(x+\delta)^p} dx
     \left(\int_0^{B\delta} V_i^{-q}(y)dy \right)^{p-1},
\end{align*}
which is bounded by a constant by \eqref{eq:3.8}. The second term is estimated
using Lemma \ref{lem:3.3} with $f(t) = g(y) /(t(1+t))$ and with $u$ and $v$ in
place of $U$ and $V$, where  
$$
 u^p(s) = U_i^p\left(\frac{As}{1+s}\right) \frac{1}{(1+s)^2} 
     \quad \hbox{and}\quad 
 v^p(s) = V_i^p\left(\frac{Bt}{1+t}\right) t^p (1+t)^{p-2}, 
$$
so that the term is bounded by, after changing the integral back to $y$,
\begin{align*}
& c \int_0^\infty \left|\frac{g(y)}{t(1+t)}\right|^p 
  V_i^p\left(\frac{Bt}{1+t}\right) (t(1+t))^p \frac{dt}{(1+t)^2} \\
&   = c \int_0^B |g(y)|^p V_i^p(Y)dY
    = c \int_{\bar t_{i-1}}^{t_i} |g(y)|^p V_i^p(|y-t_i|)dy.
\end{align*}
The condition of Lemma \ref{lem:3.3} is verified similarly; it reduces to the
condition \eqref{eq:3.9}. Putting these estimates together gives the stated
inequality \eqref{eq:3.7}. 
\end{proof}

For $U$ and $V$ being the classical GJ weight functions, the inequality 
\eqref{eq:3.7} was proved in \cite{X} under the conditions
$$
  U^p \in L^1, \qquad V^{-q} \in L^1, \qquad U(x) \le c V(x).
$$
A more general result along this line is the following:

\begin{prop} \label{prop:3.5}
Let $1 < p < \infty$. Let $U$ and $V$ be $GJ\log$ weight functions in 
Definition \ref{def:log}. Then the inequality \eqref{eq:3.7} holds if 
$p \Gamma_i(U) > -1$ or $-q \Gamma_i(V) > -1$ and $U(t) \le c V(t)$.
\end{prop}

\begin{proof}
The condition $U(x)\le cV(x)$ is equivalent to $\Gamma_i(U)\ge \Gamma_i(V)$, 
and $\gamma_i(U)\le \gamma_i(V)$ when $\Gamma_i(U) = \Gamma_i(V)$. 
That $U^p \in L^1$ implies either $p \Gamma_i(U) > -1$ or $p \Gamma_i(U)= -1$
and $p \gamma_i(U)> -1$. Hence, $p \Gamma_i(U) > -1$ or $-q \Gamma_i(V) > -1$
implies that $U^p \in L^1$ or $V^{-q} \in L^1$, respectively. We show that 
the $U$ and $V$ satisfy \eqref{eq:3.8} and \eqref{eq:3.9} under the given
conditions. 

As it is shown in \cite[p. 328]{MV}, that $p \Gamma_i(U) > -1$ and 
$-q \Gamma_i(V) > -1$ shows $U_i$ and $V_i$ satisfy 
\begin{equation} \label{eq:3.10}
\int_0^\delta U_i^p(t) dt = \CO\left(\delta U_i^p(\delta)\right) 
\quad \mbox{and} \quad
\int_0^\delta V_i^{-q}(t) dt = \CO\left(\delta V_i^{-q}(\delta)\right).
\end{equation}
Hence, $U_i(x) \le c V_i(x)$ shows that 
\begin{align*}
  \int_0^\delta \frac{U_i^p(t)}{(t+\delta)^p} dt 
      \left(\int_{0}^\delta V_i^{-q}(t) dt\right )^{p-1} 
 & \le \delta^{-p} \int_0^\delta U_i^p(U) dt \,
      \left(\delta V_i^{-q}(\delta)\right)^{p-1} \\ 
 & \le  c \,U_i^p(\delta) V_i^{-p}(\delta) \le c.
\end{align*}
Note that the part of $(\int_0^\delta V_i^{-q}(t)dt)^{p-1}$ in \eqref{eq:3.8} 
is finite since $V^{-q} \in L^1$, and it is always bounded by 
$V_i^{-p}(\delta)\delta^{p-1} = \CO(\delta^\varepsilon)$ for some 
$\varepsilon >0$. We estimate the integral of $U_i^p$ on $(\delta,1)$, 
$$
\int_\delta^1 \frac{U_i^p(t)}{(\delta+t)^p} dt \le 
        \int_\delta^1 \frac{U_i^p(t)}{t^p} dt.
$$
If $(\Gamma_i(U) - 1) p > -1$, then $U^p/t^p \in L^1$ so that \eqref{eq:3.8}
holds trivially. If $(\Gamma_i(U) - 1) p =-1$, then 
$$
 \int_\delta^1 \frac{U_i^p(t)}{t^p} dt = \int_{\delta}^1 
   \left(\log\frac{e}{t} \right)^{\gamma_i(U)p} \frac{dt}{t} = 
   \frac{\left(\log(e/\delta) \right)^{\gamma_i(U)p +1} -1}{\gamma_i(U)p+1}.
$$
In this case the fact that $V_i^{-p}(\delta)\delta^{p-1}= 
\CO(\delta^\varepsilon)$ shows that \eqref{eq:3.8} holds. Finally, if 
$(\Gamma_i(U)-1)p < -1$, we show that $\int_\delta^1 (U_i^p(t) /t^p) dt = 
\CO\left(\delta^{-p+1} U_i^p(\delta)\right)$ to finish the proof
(see also \cite[p. 328]{MV}). In this case, if $\gamma_i(U)p
\ge 0$, then since $\log (e/t)$ is decreasing,  
\begin{align*}
 \int_\delta^1 \frac{U_i^p(t)}{t^p} dt & = \int_{\delta}^1 
  t^{p \Gamma_i(U) -p} \left(\log\frac{e}{t} \right)^{p\gamma_i(U)} dt  \\
& \le c \left(\log\frac{e}{\delta} \right)^{p\gamma_i(U)}
  \int_{\delta}^1  t^{p \Gamma_i(U)-p} dt 
 \le c \,\delta^{-p+1} U_i^p(\delta),
\end{align*}
since $p \Gamma_i(U)-p +1 < 0$. If $\gamma_i(U)p < 0$, we choose $\varepsilon
> 0$ such that $(\Gamma_i(U)-1)p - \varepsilon \gamma_i(U)p <-1$. Then since
$t^\varepsilon \log(e/t)$ is increasing for $t$ close to zero, we have 
\begin{align*}
 \int_\delta^1 \frac{U_i^p(t)}{t^p} dt & = \int_{\delta}^1 
  t^{p \Gamma_i(U) -p - \varepsilon \gamma_i(U) p}
    \left(t^{\varepsilon} \log\frac{e}{t} \right)^{p\gamma_i(U)} dt  \\
 & \le c \left(\delta^{\varepsilon} \log\frac{e}{\delta} \right)^{p\gamma_i(U)}
  \int_{\delta}^1  t^{p \Gamma_i(U)-p - \varepsilon \gamma_i(U) p} dt 
  \le c \,\delta^{-p+1} U_i^p(\delta).
\end{align*}
Putting these estimates together, we have verified \eqref{eq:3.8}. The 
inequality \eqref{eq:3.9} can be verified similarly.  
\end{proof}

\medskip\noindent
{\it Remark 3.1}.
It should be pointed out that the condition $p \Gamma_i (U) > -1$ and 
$-q\Gamma_i(V) > -1$ cannot be replaced by $U^p \in L^1$ and $V^{-q}\in L^1$
without further restriction on $\gamma_i(U)$ and $\gamma_i(V)$. Assume, for 
example, that $-q\Gamma_i(V) = -1$. Then $V^{-q} \in L^1$ holds if $-q 
\gamma_i(V) < -1$, and we have 
$$
 \int_0^\delta V_i^{-q}(t) dt = \int_0^\delta 
  \left(\log\frac{e}{t} \right)^{-q\gamma_i(V)} \frac{dt}{t} = 
  \frac{\left(\log(e/\delta) \right)^{-q \gamma_i(V)+1}}{q\gamma_i(V)-1}.
$$
Here $\Gamma_i(U) \ge \Gamma_i(V) = 1/q$. If $\Gamma_i(U) = 1/q$ then 
$(\Gamma_i(U)-1)p = -1$ and the left hand side of \eqref{eq:3.8} is bounded by
$$
 \frac{\left(\log(e/\delta) \right)^{\gamma_i(U)p +1} -1}{\gamma_i(U)p+1}
\left(\frac{\left(\log(e/\delta) \right)^{-q \gamma_i(V)+1}}
{q\gamma_i(V)-1}\right)^{p-1}
 = c \left(\log\frac{e}{\delta} \right)^{p(\gamma_i(U)-\gamma_i(V)+1)},
$$
which is bounded only if $\gamma_i(U)-\gamma_i(V)+1 \le 0$. \qed
\medskip

Using the same argument, by induction if necessary, one can establish the
similar result for weight functions of the type $t^{\Gamma_i} 
(\log_k(e/t))^{\gamma_i}$, where $\log_k t = \log \log \ldots \log t$ 
($k$ fold of $\log$). To state a more general result, we need the following
definition: 

\begin{defn} \label{def3.6}
A function $S$ is called slowly varying if for any $\varepsilon > 0$, 
$t^\varepsilon S(t) \to \infty$ and $t^{-\varepsilon} S(t)\to 0$ as $t \to 
\infty$.
\end{defn}

For example, $S(t) = \log(e/t)$ is a slowly varying function and so is a power
of $S(t)$. Also, $(\log \log e/t)^\gamma$, and more generally, the powers of 
$\log_k e/t$ are all slowly varying functions. A slightmodification of the
proof of Proposition \ref{prop:3.5} gives the following result: 

\begin{prop} \label{prop:3.7}
Let $S_{U_i}$ and $S_{V_i}$ be slowly varying functions such that either 
they are increasing functions or, for any $\varepsilon >0$, $t^{-\varepsilon} 
S_{U_i}(1/t)$ and $t^{-\varepsilon} S_{V_i}(1/t)$ are increasing for small 
$t>0$. Then the inequality \eqref{eq:3.7} holds if $U$ and $V$ are GJ2 weight
with
$$
U_i(t) = t^{\Gamma_i(U)} S_{U_i}(1/t) 
\quad \hbox{and} \quad  
V_i(t) = t^{\Gamma_i(U)} S_{V_i}(1/t), \qquad 0 \le i \le r+1, 
$$
such that $p \Gamma_i(U) > -1$, $-q \Gamma_i(V) > -1$, and $U(x) \le c V(x)$.
\end{prop}

The assumption that $U$ and $V$ are in GJ2 replaces the condition 
\eqref{eq:3.10}. The Proposition \ref{prop:3.5} corresponds to the case that 
both $S_{U_i}$ and $S_{V_i}$ are of the form  $(\log (e t))^{p\gamma_i}$. If 
$\gamma_i \le 0$, this function is increasing; if $\gamma_i \ge 0$, then 
$\varepsilon^{-\varepsilon} (\log (e /t))^{p\gamma_i}$ is increasing.

Let us mention that the definition of the slowly varying functions may be 
different in the literature. For example, in \cite{Se}, it is defined as the 
functions that satisfy the relation
$$
  \lim_{t \to \infty} \frac{S(\lambda t)}{S(t)} = 1, \qquad \hbox{for each} 
   \quad \lambda > 0.
$$
Clearly, our condition is more relaxed. A slowly varying function $S$ that 
satisfy the above limiting condition also satisfies the following property 
(\cite{Se}): If $S$ is defined on $[a,\infty)$, $a > 0$, then there exists a 
$b \ge a$ such that for all $t \ge b$
$$
  S(t) = \exp \left [ \eta(t) + \int_b^t \frac{\epsilon(x)}{x} dx \right ], 
$$
where $\eta$ is a bounded measurable function on $[b,\infty)$ such that
$\eta(t) \to c$ ($|c| \to \infty$), and $\epsilon$ is a continuous function
on $[b,\infty)$ such that $\epsilon(t) \to 0$ as $t \to \infty$.

\section{convergence of orthogonal series} 
\setcounter{equation}{0}

Let $d\alpha$ be a GJ measure and we assume that $\alpha'$ is in the 
form of \eqref{eq:2.3}
$$
  \alpha'(x) = h(t) \prod_{i=0}^{r+1} \alpha_i(|t-t_i|), 
       \qquad \alpha_i(t) = [\tau_i(t)]^{a_i}.
$$
Let $S_n(d\alpha;f)$ be the partial sum of the Fourier orthogonal series. 
By \eqref{eq:1.3},  
$$
S_n (d\alpha;f,x) = \sum_{k=0}^{n-1} c_k (f) p_k (d\alpha)
    = \int_{-1}^1 f(y) K_n(d\alpha;x,y) d\alpha(y). 
$$
The kernel $K_n(d\alpha;x,t)$, by the Christoffel-Darboux formula, satisfies
the formula 
$$
K_n(d\alpha;x,y) = a_{n-1}
 \frac{p_n(d\alpha;x)p_{n-1}(d\alpha;y)-p_n(d\alpha;y)p_{n-1}(d\alpha;x)}{x-y},
$$
where $a_{n-1}$ is a proper constant. Our main results on the mean convergence 
of the generalized Jacobi series are the following (cf. \cite[I, p. 246]{X}). 

\begin{thm} \label{thm:4.1}
Let $d\alpha$ be an admissible GJ measure and assume that $\alpha_i$ are
nondecreasing for $1 \le i \le r$ and that $\alpha_i\varphi^3$ are 
nondecreasing for $i = 0$ and $r+1$. Let $u, w \in GJ$. Define $U$ and $V$ by 
$$
  U^p:= w^p(\alpha' \varphi)^{-p/2} \alpha' \quad \hbox{and}
    \quad V^{-q}: = \varphi^{q} u^{-q}(\alpha'\varphi)^{-q/2} \alpha'.
$$
Assume that $U^p$ and $V^{-q}$ satisfy \eqref{eq:3.8} and \eqref{eq:3.9}. Let 
$1<p<+\infty$. Then 
\begin{equation} \label{eq:4.1}
\| S_n(d\alpha,f)w\|_{d\alpha,p} \le c\| fu\|_{d\alpha,p}
\end{equation} 
for every $f$ such that $\| fu\|_{d\alpha,p} < +\infty$ if and only if
\begin{align} \label{eq:4.2}
\begin{split}
& w^p\alpha' \in L^1\quad ,\quad u^{-q} \alpha' \in L^1\ ,\\
& w^p\bigl( \alpha' \varphi \bigr)^{-p/2} \alpha'\in L^1\quad
,\quad u^{-q} \bigl( \alpha' \varphi\bigr)^{-q/2}\alpha'\in L^1
\end{split}
\end{align}
and
\begin{equation} \label{eq:4.3}
w(x) \le cu(x).
\end{equation} 
\end{thm}

\begin{proof}
Let $q_n(x)$ denote the orthonormal polynomials associated with the
measure $(1-x^2)\,d\alpha (x)$; that is, $q_n (x) = p_n (\varphi^2)\,
 d\alpha,x)$. Let
\begin{align*}
h_1(x,y) & = p_n(d\alpha,x) p_n(d\alpha,y),\\
h_2(x,y) & = \frac{F_n(x,y)}{x-y}, \qquad F_n(x,y)= 
   (1-y^2) p_n(d\alpha,x) q_{n-1} (y)\\
h_3(x,y) & = h_2(y,x).
\end{align*}
Following Pollard \cite{P}, the kernel $K_n(d\alpha;x,t)$ can be written as
$$
K_n (d\alpha; x,y) = \alpha_n h_1 (x,y) + \beta_n h_2(x,y) + \beta_n h_3 (x,y),
$$
where the numbers $\alpha_n$ and $\beta_n$ depend on $d\alpha$ and $n$. Since
$\alpha' >0$, a.e., it follows from \cite{R} that $|\alpha_n|$ and $|\beta_n|$ 
are bounded by constant independent of $n$ (cf. \cite[p. 358-360]{P}). Thus,
it is sufficient to prove that 
\begin{equation} \label{eq:4.4}
\int_{-1}^1 \Big| \int_{-1}^1 h_k (x,y) f(y)\, d\alpha (y)\Big|^p w^p(x)
\, d\alpha (x) \le c\| fu\|_{d\alpha,p}^p
\end{equation}
for $k=1,2$ and $3$ under the conditions \eqref{eq:4.2} and \eqref{eq:4.3}.

The bound of $p_n(d\alpha)$ in Lemma \ref{lem:2.4} shows, in particular, that
$$
  |p_n(d\alpha;x)| \le c \left[1 + (\alpha'(x) \varphi(x))^{-1/2} \right]. 
$$
A similar estimate also applies to $q_n(y)$. Applying the H\"older 
inequality to the inner integral and then using the bounds of orthogonal 
polynomials, it follows readily that the inequality \eqref{eq:4.4} holds
for $k=1$ under the condition \eqref{eq:4.2}. 

To prove \eqref{eq:4.4} for $k=2$, we use Theorem \ref{thm:3.4}. First, 
by Lemma \ref{lem:2.5} (recall $\alpha' \in GJ2$) and the fact that 
$S_n(d\alpha;f,x)$ is a polynomial of degree $n$, it is sufficient to prove 
that 
\begin{equation}\label{eq:4.5}
\int_{-1}^1 \Big|\int_{-1}^1 f(y)\frac{F_n(x,y)}{x-y}d\alpha(y)\Big|^p 
  w^p(x)
   \chi_{\Delta_n(\varepsilon)}(x)\, 
 d\alpha (x) \le c\| fu\|_{d\alpha,p}^p 
\end{equation}
under the conditions \eqref{eq:4.2} and \eqref{eq:4.3}. Lemma \ref{lem:2.4} 
shows that 
$$
|p_n(x)| \le c \bigl( \alpha'(x) \varphi(x)\bigr)^{-1/2}, \qquad 
   x \in \Delta_n (\varepsilon).
$$
Furthermore, the assumption that $\alpha_i$ is nondecreasing for $0 \le i \le 
r$ and $\alpha_i \varphi^3$ is nondecreasing for $i = 0$ and $i=r+1$ shows 
that $q_n$ is bounded by
\begin{equation}\label{eq:4.6}
 |q_n(y)| \le c \bigl( \alpha'(n,y) \varphi^3(n,y)\bigr)^{-1/2} \le c 
  \bigl( \alpha'(y) \varphi^3(y)\bigr)^{-1/2}, 
  \qquad y \in (-1,1). 
\end{equation}
Hence, by the definition of $F_n(x,y)$, the left hand side of \eqref{eq:4.5} 
is bounded by
$$
c \int_{-1}^1  \left|\int_{-1}^1 
  \frac{f(y)(\alpha'(y)\varphi^3(y))^{-1/2}\phi_n(y)}{x-y} \varphi^2(y)
 d \alpha(y)
 \right|^p (\alpha'(x) \varphi(x))^{-p/2} w^p(x) d\alpha (x), 
$$
where $\phi_n(y)$ is a function bounded by a constant independent of $n$,
which is bounded by $c \|fu\|_{d\alpha,p}^p$ upon using Theorem \ref{thm:3.4}
with $f(\alpha'\varphi^{3})^{-1/2} \phi_n$  in place of $g$. 

For $k=3$, we use a dual argument and derive the desired bound from the
case $k=2$. This argument does not depend on the fact that
our weight functions are the generalized Jacobi ones. We refer to 
\cite[p. 889]{X0} for the details. Thus, the proof for the sufficient 
part is completed. 

The conditions \eqref{eq:4.2} are proved to be necessary for general weight 
functions in \cite{MNT}. The condition \eqref{eq:4.3} is necessary also for
general weight functions as the proof in \cite[p. 250]{X} shows. 
\end{proof}

Since $S_n(d\alpha,f)$ is a projection operator, Theorem \ref{thm:4.1} and 
Weierstrauss theorem give the following corollary:

\begin{cor}
Under the assumptions of Theorem \ref{thm:4.1} 
$$
\lim_{n\to\infty} \| \bigl( S_n(d\alpha,f)-f\bigr) w\|_{d\alpha,p} =0
$$
for every $f$ such that $\| fu\|_{d\alpha,p} < \infty$ if and only if
\eqref{eq:4.2} and \eqref{eq:4.3} hold.
\end{cor}

\medskip\noindent
{\it Remark 4.1}.
The assumption that $U^p$ and $V^{-q}$ satisfy \eqref{eq:3.8} and 
\eqref{eq:3.9} already implies that $U^p \in L^1$ and $V^{-q} \in L^1$. 
Hence, the condition $w^p(\alpha'\varphi)^{-p/2} \alpha'\in L^1$ in 
\eqref{eq:4.2} is redundant for the sufficient part. We include it since 
it is also a necessary condition. Note also that $u^p(\alpha'\varphi)^{-p/2} 
\alpha' = \varphi^{-q} V^{-q}$ so that the condition 
$u^p(\alpha'\varphi)^{-p/2} \alpha' \in L^1$ in \eqref{eq:4.2} is still 
needed. 

\medskip\noindent
{\it Remark 4.2}.
Since $\alpha_i(x) = [\tau_i(x)]^{a_i}$ and $\tau_i$ are nondecreasing by
definition, the assumption on the nondecreasing of $\alpha_i$ holds if $a_i$ 
are nonnegative. This assumption is not needed for the ordinary GJ weight
functions in \eqref{eq:1.1}. See the discussion after the proof of 
Corollary \ref{cor:4.2}.
\medskip

\begin{cor} \label{cor:4.2}
Let $d\alpha$ be an admissible GJ measure. Assume that $\alpha'$, $u$ and 
$v$ are in $GJ\log$ such that, for $1 \le i \le r$, $\Gamma_i(\alpha') > 0$ 
or $\Gamma_i(\alpha') = 0$ and $\gamma_i(\alpha) \le 0$. Let $U$ and $V$ be 
defined as in Theorem \ref{thm:4.1} and assume that $p\Gamma_i(U) > -1$ and 
$-q \Gamma_i(V)> -1$ for $0 \le i \le r+1$. Then the inequality \eqref{eq:4.1}
holds if and only if \eqref{eq:4.2} and \eqref{eq:4.3} hold.
\end{cor} 

\begin{proof}
The assumption on $\alpha'$ means that $\alpha'_i(t) = t^{\Gamma_i(\alpha')} 
\left(\log (e/t)\right)^{\gamma_i(\alpha')}$. It is easy to see that this 
function is nondecreasing if $\Gamma_i(\alpha') > 0$ or if $\Gamma_i(\alpha')
=0$ and $\gamma_i(\alpha')\le 0$. Since $\alpha' \in L^1$ implies that 
$\Gamma_i(\alpha') > -1$, this shows that $\alpha_i' \varphi^3$ is increasing
for $i = 0$ and $r+1$ without further conditions. The assumption that 
$p\Gamma_i(U) > -1$ and $-q \Gamma_i(V) > -1$ implies, by Proposition 
\ref{prop:3.5}, that $U^p$ and $V^{-q}$ satisfy \eqref{eq:3.8} and 
\eqref{eq:3.9}. 
\end{proof}

Again, we note that the $p\Gamma_i(U) > -1$ implies that $U^p \in L^1$ so
that part of the \eqref{eq:4.2} is redundant. If all $\gamma_i =0$ in the 
above corollary, then it deals with the ordinary GJ weight function in 
\eqref{eq:1.1}. In that case, the assumption $p\Gamma_i(U) > -1$ and 
$-q \Gamma_i(V)> -1$ are not needed, and the result was proved in \cite{Ba}
for the case of $w=u$ and in \cite{X} for the case of $w \ne u$. However, 
the results there were proved without the assumption that $\Gamma_i \ge 0$ 
for $1 \le i \le r$ (recall that $\alpha'\in L^1$ implies $\Gamma_i > -1$). 
The additional assumption is used to ensure that \eqref{eq:4.6} holds for 
all $y \in (-1,1)$. For weight functions in the Corollary \ref{cor:4.2}, it 
is possible to follow the proof in \cite{X} to remove the additional 
assumption. This amounts to apply the current proof to $f(y) 
\chi_{\tau_{n}(\varepsilon)}$, where $\tau_n(\varepsilon)$ is the 
set  
$$
  \tau_n(\varepsilon) := [-1,1] \setminus \bigcup_{i \in \sigma}
      [t_i - \varepsilon n^{-1}, t_i +\varepsilon n^{-1}],  
$$
in which $\sigma = \{i: \Gamma_i(\alpha') <0 \quad \hbox{or} \quad  
 \Gamma_i(\alpha') =0, \gamma_i(\alpha') >0\}$ (since the inequality 
\eqref{eq:4.6} holds for $y \in \tau_n(\varepsilon)$), and show that 
\begin{align*}
\int_{-1}^1 \Big| \int_{-1}^1 h_2 (x,y)f(y) \bigl( 1-\chi_{\tau_n (d)}(y)
  \bigr) \, d\alpha (y)\Big|^p w^p(x) \chi_{\Delta_n (\varepsilon)} (x) \, 
 d\alpha(x)\le c\| fu\|_{d\alpha,p}^p 
\end{align*}
separately. The second part involves difficult estimates that have to be 
worked out. It is not clear how to extend this part to general GJ weight 
functions in Theorem \ref{thm:4.1}.  

\medskip\noindent
{\it Remark 4.3}. Under the additional assumption that $f \alpha_i$ is 
continuous locally at $x = t_i$, $1 \le i \le r$, we can remove the condition 
that $\alpha_i$ are nondecreasing for $1 \le i \le r$ from the assumption of 
the Theorem 4.1. We then have $\| S_n(d\alpha,f)w\|_{d\alpha,p}$ is uniformly 
bounded for every $f$ such that $\| fu\|_{d\alpha,p} < +\infty$ and $f$ locally
continuous at $t_i$ under the conditions (4.2) and (4.3). For the proof, we 
can assume that $\alpha_i$ is decreasing, since the other case has been 
settled in the proof of Theorem 4.1. Then it follows from the estimate in 
Lemma 2.5 that $|(1-y^2)q_{n-1}(y)| \le c$. Hence, we can use the fact that 
$$
 \lim_{n \to \infty} n \int_{t_i - n^{-1}}^{t_i+n^{-1}} |f(t) \alpha'(t)|dt 
    = |f(t_i) \alpha'(t_i)|,
$$
which holds since $f\alpha'$ is locally continuous at $t_i$ to deal with
the integral over $[-1.1]\setminus \tau_n(\varepsilon)$. We omit the details.
\medskip

Let us mention, however, the following result in which $f$ is replaced by
$f\chi_{\Delta_n(\varepsilon)}$ with $\Delta_n(\varepsilon)$ is defined as
in \eqref{eq:2.5}. 

\begin{thm} \label{thm:4.3}
Let $d\alpha$ be an admissible GJ measure, $u, w \in GJ$. Assume $U$ and $V$ 
are defined as in Theorem \ref{thm:4.1} and satisfy \eqref{eq:3.8} and 
\eqref{eq:3.9}. Let $1<p<+\infty$. Then 
\begin{equation} \label{eq:4.7}
\| S_n(d\alpha,f\chi_{\Delta_n(\varepsilon)})w\|_{d\alpha,p} \le 
  c\| fu\|_{d\alpha,p} 
\end{equation} 
for every $f$ such that $\| fu\|_{d\alpha,p} < +\infty$ if 
\eqref{eq:4.2} and \eqref{eq:4.3} hold. 
\end{thm}

This result will be used in the next section to prove Marcinkiewicz-Zygmund
type inequality. Its proof follows from the remarks above since the additional
$\chi_{\Delta_n(\varepsilon)}$ allows us to use the fact that the inequality
\eqref{eq:4.6} holds for $y \in \Delta_n(\varepsilon)$, so that the proof of 
Theorem \ref{thm:4.1} can be followed through without the additional assumption
on $\alpha'$ being nondecreasing.

\section{Marcinkiewicz-Zygmund inequality} 
\setcounter{equation}{0}

Recall the inequality for the quadrature sums in Lemma \ref{lem:2.6}. 
The Marcinkiewicz-Zygmund inequality is the converse inequality (cf. 
\cite[Theorem 2.1]{X0}. 

\begin{thm} \label{thm:6.1}
Let $d\alpha$ be an admissible GJ measure, $\beta$ be a GJ measure such that
$\beta' \in GJ2$, and $u$ be a GJ weight function such that $u^{1-q}\alpha'
\in GJ2\cap GJ4$. Define $U$ and $V$ by $U^q = u^{1-q}(\alpha'\varphi)^{-q/2}
\alpha'$ and $V^{-p} = \varphi^p \left(\alpha' \varphi\right)^{-p/2} \beta'$ 
and assume that they satisfy \eqref{eq:3.8} and \eqref{eq:3.9} with $p$ and 
$q$ exchanged. Let $P\in\Pi_{n-1}$ and $1<p<+\infty$. Then 
\begin{equation} \label{eq:6.1}
\| P\|_{d\beta,p} \le c\left(\sum_{k=1}^n |P\left( x_{kn} (d\alpha)\right)|^p
  u\left(x_{kn}(d\alpha)\right) \lambda_{kn} (d\alpha)\right)^{1/p}
\end{equation}
provided
\begin{align} \label{eq:6.2}
 u\alpha' \ge c\beta', \qquad\quad
\left(\alpha' \varphi\right)^{-p/2} \beta'\in L^1. 
\end{align}
\end{thm}

\begin{proof}
We will write $x_{kn}$ for $x_{kn}(d\alpha)$ in the following. Applying 
Lemma \ref{lem:2.5} and the usual duality argument gives  
$$
 \| P\|_{d\beta,p} \le c \| P \chi_{\Delta_n(\varepsilon)}\|_{d\beta,p} 
 = c \sup_{\| g\|_{d\beta,q}=1} \int_{-1}^1 P(t)
      \chi_{\Delta_n(\varepsilon)}(t) g(t)\, d\beta. 
$$
By the orthogonality, the Gauss-Jacobi quadrature and the H\"older inequality, 
and Lemma \ref{lem:3.1}, we have 
\begin{align*}
&\int_{-1}^1 P(t) \chi_{\Delta_n(\varepsilon)}g(t)\,d\beta
 = \int_{-1}^1 P(t) S_n (d\alpha; \chi_{\Delta_n(\varepsilon)}g\beta'
   \alpha^{\prime -1},t)\, d\alpha\\
& = \sum_{k=1}^n P(x_{kn}) S_n (d\alpha;\chi_{\Delta_n(\varepsilon)}
    g\beta'\alpha^{\prime -1},x_{kn})\lambda_{kn} (d\alpha)\\
& \le c\left(\sum_{k=n}^n |P(x_{kn})|^p \lambda_{kn} (d\alpha)
u(n, x_{kn})\right)^{1/p} \| S_n (d\alpha; g\chi_{\Delta_n(\varepsilon)}
\beta'\alpha^{\prime -1}) u^{-1/p} \|_{d\alpha,q},
\end{align*}
where in the last step Lemma \ref{lem:2.6} is used with $q$ and
$u^{1-q}\alpha'$ in place of $p$ and $v$ (recall that $u^{1-q}\alpha'\in 
GJ2\cap GJ4$). We now apply Theorem \ref{thm:4.3} with 
$g \beta'\alpha^{\prime -1}$, $u^{-1/p}$, 
$(\beta^{\prime -1} \alpha')^{1/p}$, and $q$ in place of $f$, $w$, $u$
and $p$, and conclude that
$$
\| S_n (d\alpha;g\chi_{\Delta_n(\varepsilon)}\beta' \alpha^{\prime -1}) 
u^{-1/p} \|_{d\alpha,q} \le c\| g\|_{d\beta,q}.
$$
The conditions \eqref{eq:4.2} and \eqref{eq:4.3} become conditions 
\eqref{eq:6.2} under this substitution. This completes the proof.
\end{proof}

\medskip\noindent
{\it Remark 6.1}. The condition $U^q \in L^1$ is not included in the 
condition \eqref{eq:6.2}, since it is a consequence of $U^q$ and $V^{-p}$
satisfying \eqref{eq:3.8} and \eqref{eq:3.9} with $p$ and $q$ exchanged.
Also, $u^{1-q}\alpha' \in L^1$ is not included, since it is a consequence of
$u^{1-q}\alpha' \in GJ2$.

\medskip

To use this result it is necessary to choose a weight function $u$. We need
to choose it so that $u^{1-q} \alpha' \in L^1$ (implied by $u^{1-q} \alpha'
\in GJ2$), $U^q \in L^1$ and $u \alpha' \ge c \beta'$. One choice is as 
follows: For $d\alpha, d\beta$ being GJ measure, define 
\begin{equation} \label{eq:6.3}
\sigma' = \min\{\alpha', \beta', \varphi^{-1}\}
\end{equation}
and choose $u = \alpha'\sigma'^{-1}$. By the definition of GJ weight functions
in \ref{def2.1}, it is easy to see that $d \sigma$ is also a GJ measure. 
This substitution is made in the following corollary. 



\begin{cor} \label{cor:6.2}
Let $d\alpha$ be an admissible $GJ\log$ measure and $\beta'$ be a $GJ\log$ 
measure. Let $P\in\Pi_{n-1}$ and $1<p<+\infty$. Then 
\begin{equation} \label{eq:6.4}
\| P\|_{d\beta,p} \le c\left(\sum_{k=1}^n |P\left( x_{kn} (d\alpha)\right)|^p
        \lambda_{kn} (d\sigma)\right)^{1/p},
\end{equation}
if $\Gamma_i\left((\alpha'\varphi)^{-p/2}\beta'\right)>-1$ for $1 \le i \le r$ 
and $\left(\alpha' \varphi\right)^{-p/2} \beta'\in L^1$. 
\end{cor} 

\begin{proof}
By the definition of $\sigma'$, it is evident that $\beta' \sigma'^{-1}\le c$,
$\alpha'\sigma'^{-1}\le c$ and $\varphi^{-1}\sigma'^{-1}\le c$. Consequently, 
setting $u = \alpha' \sigma^{-1}$, the condition $u\alpha' \ge c \beta'$ holds
trivially. Since $\alpha' \in GJ2$ implies that $\Gamma_i(\alpha') > -1$, the 
definition of $\sigma'$ shows that $\Gamma_i(\sigma') > -1$. Moreover, 
since $u^{1-q} \alpha' \le c \sigma'$, this shows that $\Gamma(u^{1-q} \alpha')
> -1$ so that $u^{1-q} \alpha' \in GJ2$. Let $U$ be defined as in the previous
theorem. Then $U^q \le c \sigma'$ and it follows that $\Gamma_i(U^q) > -1$. 
Furthermore, $\left(\alpha' \varphi\right)^{-p/2} 
\beta'\in L^1$ means that $\varphi^{-p} V^{-p} \in L^1$, so that 
$-p\Gamma_i(V) > -1$ for $i=0$ and $i=r+1$. Together with the assumption this
shows that $-p \Gamma_i(V) > -1$ for all $i$. Consequently, using Proposition
\ref{prop:3.5} with $p$ and $q$ exchanged finishes the proof.
\end{proof}

Our next step is to extend the Marcinkiewicz-Zygmund type inequality to 
include derivatives of $P$ in the right hand side. For this we need the 
definition of the Hermite interpolation polynomials. 

For a $m-1$ times differentiable function $f$, the Hermite interpolating 
polynomials corresponding to the distribution $d\alpha$, denoted by 
$H_{mn} (d\alpha,f)$, are defined to be the unique polynomial of degree at
most $mn-1$ satisfying 
\begin{align}\label{eq:6.5}
H_{mn}^{(j)} (d\alpha; f,x_{kn}) = f^{(j)} (x_{kn}),\qquad
0\le j\le m-1\ ,\ 1\le k\le n,
\end{align}
where $x_{kn} = x_{kn} (d\alpha)$ are zeros of $p_n(d\alpha)$. When $m=1$, 
$H_{mn} (d\alpha;f)$ are the Lagrange interpolating polynomials, we write 
$L_n(d\alpha;f) = H_{1,n} (d\alpha;f)$.

Let $d\alpha$ be an admissible GJ measure. Associated with $d\alpha$, we 
let $v$ a GJ weight function such that 
\begin{align}\label{eq:6.6}
v^{-1} (x) \le c\quad\hbox{and}\quad v^* (x):=\alpha'(x) \varphi(x)
  v^{-1} (x) \le c.
\end{align} 
Following the proof in \cite[Theorem 3.3]{X} we prove:

\begin{thm} \label{thm:6.3}
Let $m\ge1$, $P\in \Pi_{mn-1}$, and $1<p<+\infty$. Let $d\alpha$ be an 
admissible $GJ$ measure, $d\beta$ be a GJ measure such that $\beta'\in GJ2$.
Let $v$ be a GJ weight function satisfying \eqref{eq:6.6}. Let $u \in GJ$ 
such that $u^{1-q} v^{(m-1)q/2} \alpha' \in GJ2\cap GJ4$ and $u \alpha' 
\in GJ2 \cap GJ4$. For $j =1, 2 ,\ldots, m$, define $U$ and $V$ 
by $U^q = u^{1-q}v^{(j-1) q/2}(\alpha'\varphi)^{-q/2} \alpha'$
and $V^{-p} = \varphi^p \left(\alpha' \varphi\right)^{-jp/2} \beta'$ and 
assume that they satisfy \eqref{eq:3.8} and \eqref{eq:3.9} with $p$ and $q$ 
exchanged. Then 
\begin{align} \label{eq:6.7}
\| P\|_{d\beta,p} \le c\left(\sum_{j=0}^{m-1} \frac{1}{n^{jp}}
  \sum_{k=1}^n \left|\left(\varphi(x_{kn})\right)^j P^{(j)} (x_{kn})\right|^p  
 u(x_{kn}) \lambda_{kn} (d\alpha) 
  \right)^{1/p} 
\end{align}
where $x_{kn} = x_{kn} (d\alpha)$, provided
\begin{align} \label{eq:6.8}
u^{1-q} v^{(m-1)q/2}\bigl( \alpha' \varphi\bigr)^{-q/2}\alpha' \in L^1,\qquad
u\alpha' \ge c\beta' (v^*)^{-(m-1)p/2},
\end{align}
and 
\begin{align} \label{eq:6.9}
\left(\alpha' \varphi\right)^{-mp/2} \beta' \in L^1.
\end{align}
\end{thm}

\begin{proof}
We use induction. The case $m=1$ is precisely Theorem \ref{thm:6.1}. Let
us write $\eqref{eq:6.8}_m$ and $\eqref{eq:6.9}_m$ to denote the dependency 
of these conditions on $m$. We first show that $\eqref{eq:6.8}_m$ and
$\eqref{eq:6.9}_m$ imply $\eqref{eq:6.8}_{m-1}$ and $\eqref{eq:6.9}_{m-1}$. 
For \eqref{eq:6.8}, this follows as an immediate consequence of \eqref{eq:6.6}.
For \eqref{eq:6.9}, we use the fact that, for $1 \le i \le r$, if 
$\alpha_i'(x)\ge c$ then $(\alpha'\varphi)^{-(m-1)p/2} \beta' \le c \beta'$ on
$[\bar t_{i-1}, \bar t_i]$, and if $\alpha'(x) \ge c$ then 
$(\alpha'\varphi)^{-(m-1)p/2} \beta' \le (\alpha'\varphi)^{-mp/2} \beta'$ on 
$[\bar t_{i-1}, \bar t_i]$, and similar inequalities for $i = 0$ and $r+1$.

Suppose the theorem has been proved for polynomials in $\Pi_{(m-1)n-1}$ with 
$m \ge 2$ and assume that $P\in \Pi_{mn-1}$. By the interpolation property
of the Hermite interpolation 
$$
P(x) - H_{n,m-1} (d\alpha;P,x) = p_n^{m-1} (d\alpha;x) Q_n (x)
$$
where $Q_n \in \Pi_{n-1}$. Using the bound of $p_n(d\alpha;x)$ on 
$\Delta_n(\varepsilon)$ and Lemma \ref{lem:2.5},  
\begin{align*}
&\| P-H_{n,m-1} (d\alpha;P) \|_{d\beta,p}^p
 \le c\| (\alpha' \varphi)^{-(m-1)p/2} Q_n\|_{d\beta,p}^p\\
& \qquad \le c\sum_{k=1}^n |Q_n(x_{kn})|^p 
\bigl( v_n (x_{kn})\bigr)^{-(m-1)p/2}u(x_{kn}) \lambda_{kn} (d\alpha),
\end{align*}
where the last inequality follows from Theorem \ref{thm:6.1} with 
$(\alpha' \varphi)^{-(m-1)p/2}\beta'$ in place of $\beta'$ and 
$uv^{-(m-1)p/2}$ in place of $u$. The definition of $Q$ shows that 
$$
Q_n (x_{kn}) = \frac{P^{(m-1)} (x_{kn}) - H_{n,m-1}^{(m-1)} (d\alpha;P,x_{kn})}
{(m-1)! \left[ p'_n (d\alpha;x_{kn})\right]^{m-1} }.
$$
Thus, using Lemma \ref{lem:2.4}, we can estimate the sum in two terms. The 
first one is bounded by, upon using $v^*(x) \le c$, 
$$
\frac{1}{n^{(m-1)p}} \sum_{k=1}^n | \left( \varphi(x_{kn})\right)^{m-1}
  P^{(m-1)} (x_{kn})|^p u( x_{kn}) \lambda_{kn} (d\alpha),
$$
which give the $j=m-1$ term in the right hand side of \eqref{eq:6.7}. The 
second one is bounded by,
\begin{align*}
&\frac{1}{n^{(m-1)p}} \sum_{k=1}^n |\left(\varphi(x_{kn})\right)^{m-1}
H_{n,m-1}^{(m-1)} (d\alpha;P,x_{kn})|^p \\
& \qquad\qquad \times (v^*(n,x_{kn}))^{(m-1)p/2} 
u(x_{kn}) \lambda_{kn} (d\alpha),
\end{align*}
which, by the inequality in Lemma \ref{lem:2.6} with 
$u\varphi^{m-1}(v^*)^{(m-1)p/2} \alpha'$ in place of $v$, which is in 
$GJ4\cap GJ2$ since $u \alpha'$ is, is bounded by 
\begin{align*} 
\frac{c}{n^{(m-1)p}} \| \varphi^{m-1}
H_{n,m-1}^{(m-1)} (d\alpha;P) u^{1/p} \|_{d\alpha,p}^p
\le c\| H_{n,m-1} (d\alpha;P)u^{1/p} \|_{d\alpha,p}^p,
\end{align*}
where the second inequality follows from the Bernstein-Markov inequality
in Lemma \ref{lem:2.7}. Since $H_{n,m-1} (d\alpha;P) \in \Pi_{(m-1)n-1}$, 
by induction with $(v^*)^{(m-1)p/2} u \alpha'$ in place of $\beta'$, 
this term is bounded by the right hand side of \eqref{eq:6.7} with $m-1$ 
replaced by $m-2$. The conditions $\eqref{eq:6.8}_{m-1}$ and 
$\eqref{eq:6.9}_{m-1}$ under this substitution are implied by 
$u\alpha'\in L^1$, $v^*(x) \le c$ and $\eqref{eq:6.8}_m$. Thus, we have 
proved that $\| P-H_{n,m-1}(d\alpha;P)\|_{d\beta,p}$ is bounded by the right 
hand side of \eqref{eq:6.7}. Triangle inequality
$$
\| P\|_{d\beta,p} \le \| P-H_{n,m-1} (d\alpha;P)\|_{d\beta,p} +
\| H_{n,m-1} (d\alpha;P)\|_{d\beta,p}
$$ 
and induction completes the proof.
\end{proof}

The result in the theorem is given in its general form. We can choose $v$ 
and $u$ so that the conditions become easier to check. To start with, we 
choose $v$ as 
$$
   v(x) = \max\{c, \alpha'(x) \phi(x)\},  \quad \hbox{$c$ is a constant} 
$$
($c = \|h\|_\infty$), which clearly satisfies \eqref{eq:6.6}. If $\alpha' 
\in GJ$ then $v\in GJ$. Next we define 
$$
  \sigma' = \max \{\alpha' v^{(m-1)/2}, \varphi v^{(m-1)/2}, 
      (v^*)^{-(m-1)/2} \beta'\}.  
$$ 
Then $\sigma' \in GJ$. We require that $d\sigma$ is a GJ measure, that is, 
$\sigma' \in L^1$. Since $v^*(x) \le c$, $\beta' \in L^1$ shows that this
requirement put restriction on $\alpha'$. For example, if $\alpha'\in GJ\log$,
then $\sigma'\in L^1$ if 
\begin{align} \label{eq:6.10}
 \Gamma(\alpha') > \frac{-2}{m+1}, \quad 1 \le i \le r, \quad \hbox{and}\quad
 \Gamma(\alpha') > -\frac12 -\frac{1}{m+1},  \quad i=0, r+1.
\end{align}
Note that the above restriction become $\Gamma(\alpha')>-1$ if $m=1$. 
Furthermore, under the substitution $u = \alpha'^{-1} \sigma'$, the conditions
\eqref{eq:6.8} hold trivially. Indeed, the definition shows that  
$v^{(m-1)/2} \varphi^{-1}\sigma^{\prime -1}\le c$, $\alpha' 
v^{(m-1)/2} \sigma^{\prime -1} \le c$ and $\beta' (v^*)^{-(m-1)p/2} \le c
\sigma'$. The last one is precisely the second condition in \eqref{eq:6.8}
with $u= \alpha^{\prime -1} \sigma'$. Moreover
\begin{align*}
&u^{1-q} v^{(m-1)q/2} \left( \alpha' \varphi\right)^{-q/2}\alpha'\\
&\qquad = \bigl( \alpha' v^{(m-1)/2}\sigma^{\prime -1}\bigr)^{q/2} 
 \bigl( v^{(m-1)/2} \varphi^{-1} \sigma^{\prime -1}\bigr)^{q/2}
 \sigma' \le c\sigma',
\end{align*}
which shows that the first condition of \eqref{eq:6.8} holds. Furthermore,
we also have 
$$
u^{1-q} v^{(m-1)q/2} \alpha' = \bigl( \alpha' v^{(m-1)/2} \sigma^{\prime-1}
\bigr)^q \sigma' \le c\sigma',
$$
which shows that $u^{1-q} v^{(m-1)q/2} \alpha'\in L^1$. Finally, the 
substitution $u = \alpha'^{-1} \sigma'$ leads to the inequality
\begin{align} \label{eq:6.11}
\| P\|_{d\beta,p} \le c\left(\sum_{j=0}^{m-1} \frac{1}{n^{jp}}
  \sum_{k=1}^n \left|\left(\varphi(x_{kn})\right)^j P^{(j)} (x_{kn})\right|^p  
  \lambda_{kn} (d\sigma)  \right)^{1/p} 
\end{align}
where $x_{kn} = x_{kn} (d\alpha)$. 

Hence, with these choices of $u$ and $v$ we can simplify the conditions in the
previous theorem. This shows that \eqref{eq:6.11} holds essentially under the
condition $(\alpha'\varphi)^{-mp/2} \beta' \in L^1$ or a slightly stronger 
one. To make the conditions precise will require stating assumptions on $U$ 
and $V$ precisely, which can be rather involved. Instead of trying to state 
a general result, we restrict again to the $GJ\log$ case. 

\begin{thm} \label{thm:6.4}
Let $m\ge1$, $P\in \Pi_{mn-1}$, and $1<p<+\infty$. Let $d\alpha$ be an 
admissible $GJ\log$ measure, $d\beta$ be a $GJ\log$ measure such that
$\Gamma_i(\beta') > -1$, $1 \le i \le r$. Let $\sigma'$ be defined as above 
and assume \eqref{eq:6.10} so that $\sigma'\in L^1$. Then the inequality
\eqref{eq:6.11} holds provided 
$\Gamma_i\left(\alpha'^{-mp/2}\beta'\right)>-1$ for $1 \le i \le r$ 
and  $(\alpha'\varphi)^{-mp/2} \beta' \in L^1$.
\end{thm}

\begin{proof}
We take $u = \alpha'^{-1} \sigma'$ in Theorem \ref{thm:6.3} as indicated 
above. Then the condition \eqref{eq:6.8} is already satisfied. The assumption 
\eqref{eq:6.10} and $\Gamma_i(\beta) > -1$ shows that $\Gamma(\sigma') > -1$ 
for $0 \le i \le r+1$. This implies, in particular, that $u \alpha' = \sigma'
 \in GJ2 \cap GJ4$ and $u^{1-q} v^{(m-1)q/2} \alpha' \in GJ2 \cap GJ4$ since 
the latter is bounded by $c\sigma'$.  For $j =1,2,\ldots,m$, let $U$ and $V$
be defined as in Theorem \ref{thm:6.3} with $u = \alpha^{\prime -1} \sigma'$. 
Since $v_i(x) = c$ on $[\bar t_{i-1},\bar t_i]$ if $\alpha_i'\varphi \le c$ 
and $v(x) = \alpha'_i(x)\varphi(x)$ on $[\bar t_{i-1},\bar t_i]$ otherwise, 
it follows that $v^{(j-1)q/2}(x) \le c v^{(m-1)q/2}(x)$. Hence, it follows 
that 
$$
  U^q(x) \le u^{1-q} v^{(j-1)q/2} (\alpha' \varphi)^{-q/2}\alpha'
    \le u^{1-q} v^{(m-1)q/2} (\alpha' \varphi)^{-q/2}\alpha' \le c \sigma'
$$
as before, which shows that $\Gamma_i(U^p) > -1$ since $\Gamma(\sigma')>-1$.
Furthermore, if $\alpha'_i\varphi \ge c$ on $[\bar t_{i-1},\bar t_i]$ then 
$V^{-p} = \varphi^p (\alpha'\varphi)^{-jp/2} \beta' \le c \varphi^p \beta'$
on $[\bar t_{i-1},\bar t_i]$, which shows that $\Gamma_i(V^{-p}) > -1$ by
the assumption on $\Gamma_i(\beta')$; if $\alpha'_i\varphi \le c$ on 
$[\bar t_{i-1},\bar t_i]$ then $V^{-p} \le \varphi^p (\alpha'\varphi)^{-mp/2} 
\beta'$ on $[\bar t_{i-1},\bar t_i]$, which shows that $\Gamma_i(V^{-p}) > -1$
by the assumption on $\Gamma_i\left((\alpha'\varphi)^{-mp/2}\beta'\right)>-1$. 
By Proposition \ref{prop:3.5}, this shows that $U^q$ and $V^{-p}$ satisfy 
\eqref{eq:3.8} and \eqref{eq:3.9} with $p$ and $q$ exchanged. 
\end{proof}

We note that the condition $\Gamma_i\left(\alpha'^{-mp/2}\beta'\right)>-1$ 
for $1 \le i \le r$ is just a slightly stronger than that of 
$(\alpha'\varphi)^{-mp/2} \beta' \in L^1$ inside $(-1,1)$, which implies 
$\Gamma_i\left(\alpha'^{-mp/2}\beta'\right) \ge -1$. If $\gamma_i =0$, then
the two conditions are equivalent, and the result was proved in \cite{X}.

\section{Mean convergence of interpolating polynomials} 
\setcounter{equation}{0}

With the Marcinkiewicz type inequality established, the mean convergence of 
the corresponding interpolating polynomials follows right away. In the 
following we state the result for interpolating polynomials based on the 
zeros of orthogonal polynomials with respect to a $GJ\log$ weight function.
Let $s\ge 0$, $C^s [-1,1] = C^s$ denote the space of $s$ times continuously 
differentiable functions. We begin with the following fundamental result:

\begin{thm}
Let $d\alpha$ be an admissible $GJ\log$ measure, $d\beta$ be a $GJ\log$ 
measure such that $\Gamma_i(\beta') > -1$, $1 \le i \le r$ and \eqref{eq:6.10}
holds. Assume for $0 \le \ell \le m$ that 
$$
\Gamma_i\left(\alpha'^{-mp/2}\beta'\right) >-1 \quad \hbox{for 
$1 \le i \le r$} \quad \hbox{and} \quad  
(\alpha'\varphi)^{-mp/2} \varphi^{\ell p} \beta' \in L^1.
$$
Then for $f\in C^{m-1}$ 
$$
\| H_{nm} (d\alpha;f)\|_{d\beta,p} \le c n^\ell \sum_{j=0}^{m-1} 
\max_{1\le k\le n} | \varphi (x_{kn})^j f^{(j)} (x_{kn})| / n^j.
$$
\end{thm}

\begin{proof}
First let $1<p<+\infty$. Since for every fixed $d>0$, $n^{-1} \le \varphi(x)$
on $[-1 +dn^{-2},1-dn^{-2}]$, it follows from Lemma \ref{lem:2.5} and 
Lemma \ref{lem:2.7} that
\begin{align*}
\| H_{nm} (d\alpha;f)\|_{d\beta,p} & \le c\| H_{nm} (d\alpha;f)
\chi_{\Delta_n (\varepsilon)} \|_{d\beta,p}\\
& \le cn^\ell \| H_{nm} (d\alpha;f) \varphi^t \|_{d\beta,p}.
\end{align*}
We then apply Theorem \ref{thm:6.4} with $P = H_{nm} (d\alpha;f)$ and 
$\varphi^{\ell p}\beta'$ in place of $\beta'$. Since the assumption implies
that $\sigma' \in L^1$, it follows that
$$
\sum_{k=1}^n \lambda_{n} (d\sigma; x_{kn}) \le c\int_{-1}^1 d\sigma < +\infty.
$$
This establishes the stated inequality for $1<p<+\infty$. The case for $0 < p 
\le 1$ follows from an argument in \cite[p. 88]{X} which goes back to 
\cite[p. 886]{N}.
\end{proof}

Evidently, one could state such a result based on Theorem \ref{thm:6.1} with 
$u = \sigma' \alpha'^{-1}$ for more general weight functions. The conditions
on $U$ and $V$ make it less practical. 

In the case $\ell =0$, the above theorem shows the boundedness of the 
operator $H_{nm}(d\alpha;f)$ from $L^p(d\beta)$ to $C^{m-1}$. Using the 
Bernstein-Markov inequality, one gets also the boundedness of 
$\|H_{nm}(d\alpha;f)\|_{d\beta,p}$. For example, we have the following result.

\begin{thm}
Let $d\alpha$ be an admissible $GJ\log$ measure, $d\beta$ be a $GJ\log$ 
measure such that $\phi^{-k p}\beta' \in L^1$, $\Gamma_i(\beta') > -1$,
$1 \le i \le r$ and \eqref{eq:6.10} holds. Assume that 
$$
\Gamma_i\left(\alpha'^{-mp/2}\beta'\right) >-1 \quad \hbox{for 
$1 \le i \le r$} \quad \hbox{and} \quad  
(\alpha'\varphi)^{-mp/2} \varphi^{(m-k-1)p} \beta' \in L^1,
$$
where $0 \le k \le m-1$. Then 
$$
\lim_{n\to\infty} \| H_{nm}^{(k)} (d\alpha;f)-f^{(k)} \|_{d\beta,p} =0,
\qquad \forall  f\in C^{m-1}.
$$
\end{thm}

\begin{proof}
Since $H_{mn}(d\alpha;f)$ is a projector from $C^{m-1}$ to $\Pi_{mn-1}$,
we only need to estimate $\|H_{nm}^{(j)} (d\alpha;f-R_n)\|_{u,p}$, where 
$R_n$ is a polynomial of degree $n$ such that for $0 \le j \le m-1$, 
$$
\|f^{(j)} - R_n^{(j)}\|_\infty \le c E_n (f^{(m-1)})\big\slash n^{m-1-j},
\qquad \forall\ f\in C^{m-1},
$$
in which $E_n(f) = \inf_{P\in\Pi_n} \| f-P\|_\infty$. Using the 
Bernstein-Markov inequality shows that 
$$
\| H_{nm}^{(k)} (d\alpha; f-R_n)\|_{d\beta,p} \le c n^k \|H_{nm}
(d\alpha;f-R_n)\varphi^{-k}\|_{d\beta,p},
$$
which is bounded by $cE_n (f^{(m-1)})\to 0$ upon applying the previous 
theorem with $\ell = m-1 - k$ and $\beta'\varphi^{-jp}$ in place of $\beta'$.
\end{proof}

In particular, for $m=1$, this shows the convergence of Lagrange interpolation.

\begin{cor}
Let $d\alpha$ be an admissible $GJ\log$ measure, $d\beta$ be a $GJ\log$ 
measure such that $\Gamma_i(\beta') > -1$, $1 \le i \le r$. Assume that
$$
\Gamma_i\left(\alpha'^{-p/2}\beta'\right) >-1 \quad \hbox{for 
$1 \le i \le r$} \quad \hbox{and} \quad  
(\alpha'\varphi)^{-p/2} \beta' \in L^1.
$$
Then 
$$
\lim_{n\to\infty} \| L_n(d\alpha;f)-f \|_{d\beta,p} =0,
 \qquad \forall  f\in C.
$$
\end{cor}

The method also allows us to prove result concerning the best convergence 
order of the interpolating polynomial. 

\begin{thm}
Let $d\alpha$ be an admissible $GJ\log$ measure, $d\beta$ be a $GJ\log$ 
measure such that $\phi^{-k p}\beta' \in L^1$, $\Gamma_i(\beta') > -1$,
$1 \le i \le r$ and \eqref{eq:6.10} holds. Assume that 
$$
\Gamma_i\left(\alpha'^{-mp/2}\beta'\right) >-1 \quad \hbox{for 
$1 \le i \le r$} \quad \hbox{and} \quad  
(\alpha'\varphi)^{-mp/2} \varphi^{-k p} \beta' \in L^1,
$$
where $0 \le k \le m-1$. Then 
$$
\lim_{n\to\infty} \| H_{nm}^{(k)} (d\alpha;f)-f^{(k)} \|_{d\beta,p} 
    \le c E_n(f^{(m-1)})/ n^{m-k-1}, \qquad \forall  f\in C^{m-1}.
$$
\end{thm}

These results include many special cases considered by various authors. 
See, for example, discussions in \cite{MV,N,X}.

One can also apply the approach to other type of interpolation processes,
for example, to Hermite-Fe\'er interpolation polynomials and to truncated
Hermite interpolation polynomials. See, for example, the discussion in 
\cite[Section 4]{X}.

\enddocument